\documentclass[11pt]{article}
\usepackage{tikz-cd}
\usepackage{amsmath}
\usepackage{amsthm}
\usepackage{amsfonts}
\usepackage{amssymb}
\usepackage{color}
\usepackage{graphicx}
\usepackage{subcaption}
\usepackage{url}
\usepackage{microtype}
\usepackage{amsmath}
\usepackage{amsfonts}
\usepackage{graphicx}
\usepackage{booktabs} 

\usepackage{hyperref}

\usepackage[top=2.5cm, bottom=3.5cm, left=2.5cm, right=2.5cm]{geometry} 
\renewcommand{\epsilon}{\varepsilon}

\usepackage{mathtools}

\usepackage[noend]{algorithmic}
\usepackage{algorithm}
\floatname{algorithm}{Algorithm}

\title{Ball mapper: a shape summary for topological data analysis.}

\author{Pawe{\l} D{\l}otko \\
        Department of Mathematics, Swansea University}

\date{ \today}

\begin{document}         
\maketitle

\begin{abstract}
Topological data analysis provides a collection of tools to encapsulate and summarize the shape of data. Currently it is mainly restricted to \emph{mapper algorithm} and \emph{persistent homology}. In this paper we introduce new mapper--inspired descriptor that can be applied for exploratory data analysis.
\end{abstract}

\section{Introduction}
\label{sec:intro}
\emph{Topological Data Analysis} (TDA) is a new field in data science. It aims to use rigorous methods of topology, as \emph{Reeb} and \emph{mapper graphs} as wel as \emph{persistent homology}, to analyze complicated, large, noisy and high dimensional data. It has been successfully applied to a large number of tasks including analysis of genetic data sets~\cite{gunnar-brest-cancer}, fraud detection, biology and more. We suggest the reader to consult~\cite{gunnar_bams} as well as Gunnar Carlson's blog~\cite{gunnar_blog} for further details on usage of mapper algorithm and TDA.

Despite its success, TDA for more than ten years, provide tools restricted only to persistent homology and mapper algorithm. There has been a considerable progress both in theoretical understanding and practical computations of both, however there are no new rigorous shape descriptors.

In this paper we propose new mapper--inspired descriptors of data called a \emph{Ball mapper algorithm} (BM algorithm). It rigorously encapsulate both the local and the global structure of given data sets and can be used in exploratory data analysis. Given its simplistic construction it is easy to compute and analyze and it is scalable.

\section{Preliminaries}
\label{sec:preliminaries}
In this section we will briefly introduce the Mapper algorithm (MA) and Persistent homology (PH).

\subsection{Persistent homology}
\label{sect:PH}
In this paper we are using \emph{persistent homology} as a tool required to show some of theoretical property of the BM algorithm. As it is a standard topic discussed in many textbooks and papers, we will present only a basic introduction.

The main building block we need is a \emph{abstract simplicial complex} $\mathcal{K}$. It is a set of sets such that for every $s \in \mathcal{K}$ and every $t \subset s$, $t \in \mathcal{K}$. Elements of $\mathcal{K}$ are called \emph{simplices}. Starting from a point cloud $X$ one can obtain an abstract simplicial complex for instance by using \emph{Vietoris-Rips} (VR) complex construction. We place balls of radius $r$ centered at every $x \in X$. A simplex in VR  complex is supported in balls centered in points $x_1,\ldots,x_n \in X$ if $B(x_i,r) \cap B(x_j,r) \neq \emptyset$ for every pair of $i,j \in \{ 1,\ldots,n \}$. The VR complex of a radius $r$ supported on point cloud $X$ will be refereed to as $VR(X,r)$.

In addition one can let $r$ vary between $0$ and $+\infty$.  Clearly $VR(X,r) \subset VR(X,r')$ for $r < r'$, which gives a \emph{filtration} (nested sequence) of VR complexes. The sequence of inclusions between complexes induce a sequence of isomorphisms in homology groups $f_{r,r'} : H_p(VR(X,r)) \rightarrow H_p(VR(X,r'))$. Please consult~\cite{herbert} for a definition of homology groups. The \emph{persistent homology groups} are the images of $f_{r,r'}$, (see~\cite{herbert} for details).

Intuitively persistent homology track the evolution of the number of connected components and holes in filtered complexes (as VR complexes obtained for increasing sequence of r's).

\subsection{The Mapper algorithm}
\label{sec:mapper}
\begin{figure}[ht]
\vskip 0.2in
\begin{center}
\includegraphics[scale=1]{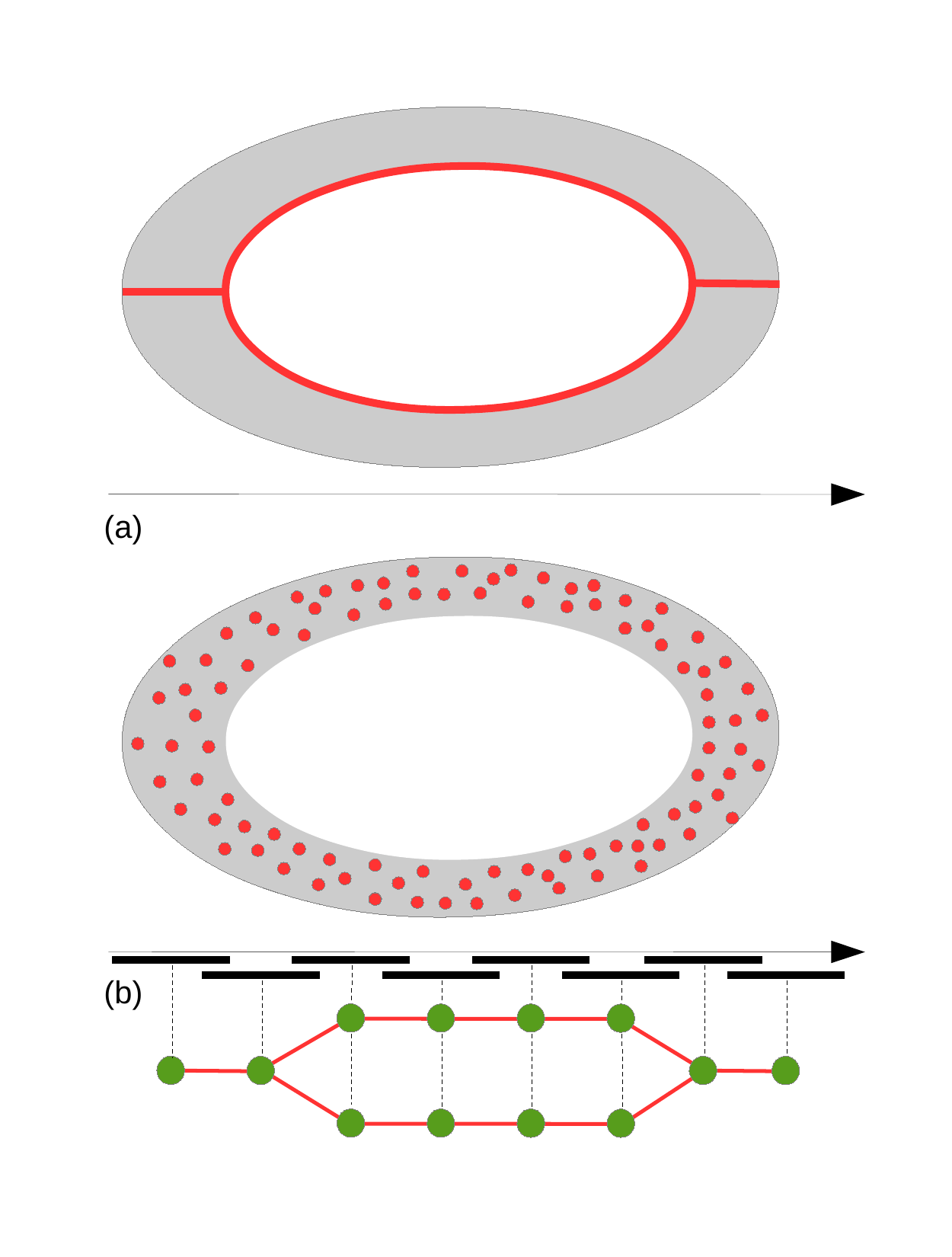}
\caption{\emph{(top, part a) Idea of a construction of a Reeb graph, (bottom, part b) idea of Mapper algorithm.}}
\label{fig:Reeb_and_mapper}
\end{center}
\vskip -0.2in
\end{figure}

The Mapper algorithm, refereed to in this paper as \emph{Conventional Mapper (CM)}, is motivated by the idea
of a \emph{Reeb graph}. Suppose we are given a manifold $\mathcal{M}$ and a continuous real value function $f :
\mathcal{M} \rightarrow \mathcal{R}$. We define an equivalence relation $\sim $ in $X$ in the following way: for every
pair of points $x,y \in \mathcal{M}$ we say that $x \sim y$ if $f(x) = f(y)$, and $x$ and $y$ are
in the same connected component of $f^{-1}(x) (=f^{-1}(y))$. The quotient space
$\mathcal{M}/_{\sim}$ is called a \emph{Reeb graph}~\cite{reeb_graph}. Please consult the upper part (a) of the
Figure~\ref{fig:Reeb_and_mapper} for the idea: the oval shape represents a manifold $\mathbb{M}$. The function $f:
\mathbb{M} \rightarrow \mathbb{R}$ is defined to be the x-coordinate of a point (equivalently, a
projection of a point to the x-axis). For every $x$, $f^{-1}(x)$ is either empty (outside the vertical stripe occupied by
$\mathbb{M}$), consist of one, (on both sides of $\mathbb{M}$), or two (in the middle), connected
components. The Reeb graph is visualized with the red lines. It is a one dimensional summary of the shape of $\mathbb{M}$.

Mapper algorithm was originally introduced in~\cite{mapper}. It is a patented core of Ayasdi~\cite{ayasdi} data analysis software. The idea of Mapper is based on consideration whether one construct a discretization of a Reeb graph given a finite set $X$ sampled iid from $\mathbb{M}$. In this case, generically $f^{-1}(x) = \emptyset$, therefore the inverse images of points are replaced by inverse images of a collection of overlapping closed intervals $I_1,\ldots,I_n$ that cover the range of $f$. As $f^{-1}(I_k) \subset X$ is a discrete space instead of connected components we consider clusters\footnote{The clustering algorithm and its parameters are parameters of mapper.} in $f^{-1}(I_k)$. The collection of clusters obtained for every interval $I_k$ corresponds to the vertices of the mapper graph. An edge between two vertices (connected components) $C_1$ and $C_2$ is placed if $C_1 \subset f^{-1}(I_k)$ and $C_2 \subset f^{-1}(I_{k+1})$ for $k \in \{1,\ldots,n-1\}$ and $C_1 \cap C_2 \neq \emptyset$\footnote{Note that this is possible because $I_{k} \cap I_{k+1} \neq \emptyset$.}.

This construction is a special case of a general \emph{Nerve complex} construction. Given a cover $C_1,\ldots,C_n$ of a set $X$, n-simplex is supported in the cover elements $C_{i_1},\ldots,C_{i_k}$ if and only if $C_{i_1} \cap \ldots \cap C_{i_k} \neq \emptyset$.

Often the vertices of mapper graph are colored by the attribute of interest: For instance, in~\cite{gunnar-brest-cancer} the authors have discovered a new type of breast cancer with $100\%$ survival rate by using the survival rate to color the obtained mapper graphs.

The obtained abstract graph will be refereed to in this paper as \emph{Conventional mapper} graph, or simply \emph{CM graph}. The idea of the construction is provided in the Figure~\ref{fig:Reeb_and_mapper}, part (b): The red points indicate the finite sampling $X$ of $\mathcal{M}$. Below, the real line (x-axis) is equipped with a collection of overlapping intervals (marked by bold black edges). Taking the inverse images via $f^{-1}$ of each of them and running clustering algorithm at those inverse images, provides the vertices of the mapper graph (in green at the bottom). Edges are added for every pair of vertices if the corresponding clusters have nonempty intersection.

It is worth mentioning that both CM and BM algorithms are not restricted to point clouds sampled from manifolds.

A great advantage (and handicap) of CM algorithm is its dependence on a number of parameters: \textbf{(1)} The input data $X$, \textbf{(2)} The function $f : X \rightarrow \mathbb{R}$, often called \emph{lens}, \textbf{(3)} The cover of $\mathbb{R}$, including the percentage of intersection of the intervals and \textbf{(4)} The clustering algorithm. Setting them in a ''correct'' way allow to get a great inside into a data and makes the CM very general tool. At the same time the art of setting them up is the main \emph{know-how} in using of CM. In the next section we present a mapper--inspired algorithm that can provide similar outputs requiring only one parameter.
\section{The ball-mapper algorithm}
The main observation we will build up on is that all the parameters of CM are used to obtain \emph{an overlapping cover} $C$ of the point cloud $X$.

In this section we will explore an alternative way of obtaining an overlapping covers of $X$. The main idea of this approach is, for a given constant $\epsilon > 0$, to find a collection of points $C \subset X$ such that the collection of balls, $B(C) = \bigcup_{x \in C} B(x,\epsilon)$ cover the entire set $X$. $B(X,\epsilon)$, refereed to as a \emph{cover vector} is collection of points in $X$ equipped with a list of ball's centers such that for every given point $x \in X$ balls of radius $\epsilon$ centered in points in $B(X,\epsilon)[x]$ cover $x$, and no other balls do.

There are a number of ways the collection of centers of balls $C$ can be selected among $X$:
\begin{enumerate}
\item By building an \emph{$\epsilon$-net on $X$} ($\epsilon$ is a parameter),
\item By using k-means clustering~\cite{k-means} or a Gaussian Mixture Model (GMM)~\cite{gmm} (number of clusters is the parameter of the method) and subsequently finding a distance $\epsilon$ of a farthest data point from the selected cluster's centers.
\end{enumerate}

As k-means and GMM's are standard tools we will skip their detailed explanation. Let us present two algorithms to construct $\epsilon$-nets: the greedy presented in Algorithm~\ref{alg:geedy_epsi_net} and max-min based presented in Algorithm~\ref{alg:max_min_epsi_net}.
While both of them produce an $\epsilon$ net, an advantage of Algorithm~\ref{alg:max_min_epsi_net} is that it can be modified so that it stops after a fixed number of iterations providing as spread collection of points as possible. Algorithm~\ref{alg:geedy_epsi_net} at the other hand is very easy to implement and distribute and its  version scales well.

\begin{algorithm}[tb]
   \caption{Greedy $\epsilon$-net}
   \label{alg:geedy_epsi_net}
\begin{algorithmic}
   \STATE {\bfseries Input:} Point cloud $X$, $\epsilon > 0$
   \STATE Mark all points in $X$ as \emph{not covered}. 
   \STATE Create initially empty cover vector $B(X,\epsilon)$.
   \REPEAT
   \STATE Pick a first point $p \in X$ that is not covered.
   \STATE For every point in $x \in B(p,\epsilon) \cap X$, add $p$ to $B(X,\epsilon)[x]$.
   \UNTIL{All elements of $X$ are covered.}
   \STATE return $B(X,\epsilon)$.
\end{algorithmic}
\end{algorithm}

\begin{algorithm}[tb]
   \caption{Max-min $\epsilon$-net}
   \label{alg:max_min_epsi_net}
\begin{algorithmic}
   \STATE {\bfseries Input:} Point cloud $X$, $\epsilon > 0$
   \STATE Pick arbitrary $p \in X$. $C = \{p\}$
   \REPEAT
   \STATE Find point $p \in X \setminus C$ that is farthest away from $C$ (if there is more than one, pick any).
   \STATE $d = dist( p,C )$
   \STATE $C = C \cup \{p\}$.
   \UNTIL{$d \leq \epsilon$}
   \STATE Create initially empty cover vector $B(X,\epsilon)$.
   \FOR {Every point $p \in C$}
        \STATE For every $x \in B(p,\epsilon) \cap X$, add $p$ to $B(X,\epsilon)[x]$.
   \ENDFOR
   \STATE Return $B(X,\epsilon)$.
\end{algorithmic}
\end{algorithm}
Given a cover vector $B(X,\epsilon)$ provided by Algorithm~\ref{alg:geedy_epsi_net} or~\ref{alg:max_min_epsi_net} one can construct a nerve or a cover $N(B(X,\epsilon))$ in the following way: The balls centers $c_1,\ldots,c_n$ corresponds to vertices of $N(B(X,\epsilon))$. A simplex $[c_{i_0},\ldots,c_{i_k}]$ is formed if there exist $x \in X$ covered by all of $c_{i_0},\ldots,c_{i_k}$\footnote{Note that $x$ is the point that belong to the intersection $B(c_{i_0},\epsilon) \cap \ldots \cap B(c_{i_k},\epsilon)$}.

The definition above naturally leads to a \emph{filtered} nerve of a cover. The filtration of a simplex $[c_{i_0},\ldots,c_{i_k}]$ being the number of points in $X$ covered simultaneously by $c_{i_0},\ldots,c_{i_k}$. Although it is a natural extension, at the moment it will not be explored here due to practical limitations in visualising edges of different weights. 

The obtained complex / graph will be refereed to as \emph{Ball mapper graph} (BM graph). The final step of a construction from $B(X,\epsilon)$ is given in the Algorithm~\ref{alg:ball_mapper_final}.

\begin{algorithm}[tb]
   \caption{Construction of Ball Mapper graph}
   \label{alg:ball_mapper_final}
\begin{algorithmic}
   \STATE {\bfseries Input:} Cover vector $B(X,\epsilon)$ from  Algorithm~\ref{alg:geedy_epsi_net} or~\ref{alg:max_min_epsi_net}.
   \STATE $V = $ cover elements in $B(X,\epsilon)$, $E = \emptyset$,
   \FOR {Every point $p \in X$}
        \STATE For every pair of cover elements $c_1 \neq c_2$ in $B(X,\epsilon)[p]$, add a (weighted) edge between vertices corresponding to the cover elements $c_1$ and $c_2$. Formally, $E = E \cup \{c_1,c_2\}$
   \ENDFOR
   \STATE Return $G = (V,E)$
\end{algorithmic}
\end{algorithm}

A simple example of the described procedure is provided in the picture below. Starting from a two dimensional window--shaped point cloud (left) a BM construction with $\epsilon = 1$ (middle left), $\epsilon = 2$ (middle right) and $\epsilon = 3$ (right) is given.

\begin{tabular}{c c c c}
      \includegraphics[width=35mm]{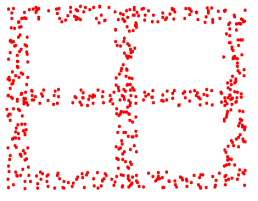} &
      \includegraphics[width=35mm]{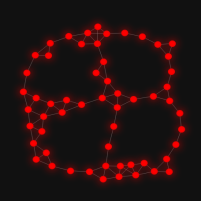} &
      \includegraphics[width=35mm]{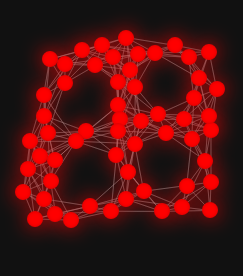} &    \includegraphics[width=35mm]{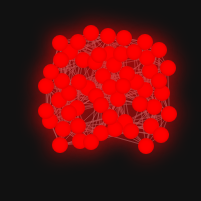} \\
\end{tabular}

Let us discuss the computational complexity of the proposed procedures. Let $B$ denotes the collection of cover ball's centers. Both for the Algorithm~\ref{alg:geedy_epsi_net} and~\ref{alg:max_min_epsi_net} for $b \in B$ we need to compute its distance to every other point in $X$. Unless an efficient spatial searching algorithms like~\cite{k_d_tree} is used, this operation requires $O( card(X)*card(B)*c )$, were $c$ is the complexity of computations of a distance between two points in $X$ (not negligible for high dimensional point clouds). We have decided not to use the algorithms as in~\cite{k_d_tree} as they are known not to provide an essential speed up in very high dimensions which are the main focus are aiming for. If they can be use, the worst case complexity presented here will be reduced.

The $card(B)$ strongly depends on the choice of $\epsilon$ parameter. For small $\epsilon$ it may be of the order of $card(X)$, which yields the worst case quadratic complexity of the two algorithms. Let us however indicate that since the aim of the presented techniques is building essentially simplified models of data, we aim for $card(B) << card(X)$, in which the complexity will substantially lower.

In the Algorithm~\ref{alg:ball_mapper_final} we first find a cover of $X$, the complexity of which is $O( card(X)*card(B)*c)$. The time complexity of the subsequent \emph{for} loop measures at the same time the complexity of the output BM graph. It can be bounded by $O(card(X) max(B(X,\epsilon)[i]))$. As $max(B(X,\epsilon)[i]) < card(B)$, $O(card(X) max(C[i]))$ is bounded by $O(card(X) card(B))$. Therefore the overall complexity of the Algorithm~\ref{alg:ball_mapper_final} is of the order $O( card(X)*card(B)*c )$.
\subsection{Multi-scale approach.}
\label{sec:mult-scale-approach}
In this section we consider the BM construction for a fixed ball's centers obtained by Algorithm~\ref{alg:geedy_epsi_net} or~\ref{alg:max_min_epsi_net} for $\epsilon_1$ and cover constructed for increasing collection of radii $\epsilon_1 \leq \epsilon_2 \leq \ldots \leq \epsilon_n$. Doing so, we will obtain a nested sequence of covers. As discussed in Section~\ref{sec:Guarantees_on_connectivity} using this technique will give theoretical guarantees on the output graph / nerve complex. The details of multi-scale BM algorithm are summarized in the Algorithm~\ref{alg:multi_scale_ball_mapper}

\begin{algorithm}[tb]
   \caption{Multi-scale Ball Mapper algorithm.}
   \label{alg:multi_scale_ball_mapper}
\begin{algorithmic}
   \STATE {\bfseries Input:} Point cloud $X$, $\epsilon_1 \leq \epsilon_2 \leq \ldots \leq \epsilon_n.$
   \STATE Pick a collection $B$ of ball centers using Algorithm~\ref{alg:geedy_epsi_net} or Algorithm~\ref{alg:max_min_epsi_net} for $X$ and $\epsilon_1$.
   \FOR {Every $\epsilon \in \{ \epsilon_1, \epsilon_2,\ldots,\epsilon_n \}$}
        \STATE For every point $p \in X$, select those elements in $B$ which are closer than $\epsilon$ to $p$ and, based on them, construct the cover vector $B(X,\epsilon)$.
        \STATE Compute graph $G_{\epsilon}$ by running Algorithm~\ref{alg:ball_mapper_final} for $B(X,\epsilon)$.
   \ENDFOR
   \STATE Return $\{ G_{\epsilon_1},\ldots, G_{\epsilon_n} \}$
\end{algorithmic}
\end{algorithm}

\subsection{Relation of CM and BM}
\label{sec:relaton_to_standard_mapper}
In this section we will discuss a relation of CM and BM algorithm.

\begin{figure}[ht]
\vskip 0.2in
\begin{center}
\includegraphics[scale=1]{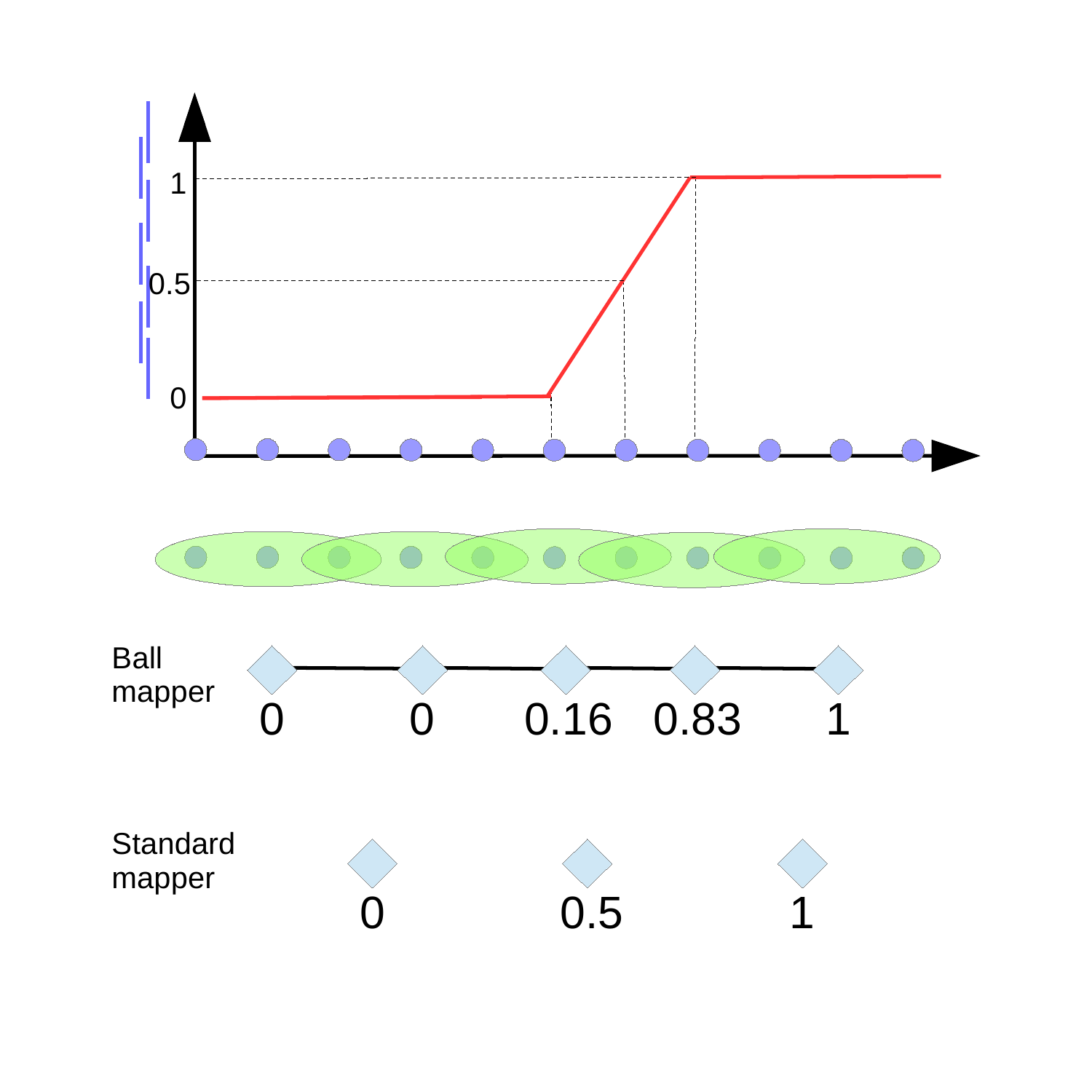}
\caption{\emph{Suppose we are given a one dimensional lens function that is initially $0$, then rapidly increase to $1$, and remain $1$ (presented in red). On the right from the $y$ axis, the cover of the range is depicted. Dots on $x$ axis represent the point cloud $X$. Note that in this example the cover is quite dense. Below the plot, with the green ovals, an example of cover used in the BM and the resulting nerve of this cover (BM graph). The numbers down to the vertices of the BM graph denote average value of function on the covered elements. Below that, a result of a CM based on the cover of the domain depicted in the left from the y-axis. It consist of three disconnected cover elements. The one having function value zero contain all the points on which function is $0$. Likewise for the points on which the function is $1$. The cover element with the value $0.5$ contain only one point, where the function value is $0.5$. The inverse images of all other intermediate cover elements are empty, what makes the CM graph disconnected.}}
\label{fig:difference_between_mappers}
\end{center}
\vskip -0.2in
\end{figure}

The key factor is the lens function $f$ used in CM. If it is not continuous, then there is no correspondence between points that are close in the domain and range of $f$. Therefore the CM and BM graphs may be arbitrarily different and nothing can be said in this case.

Let us therefore suppose that the lens function $f : X \rightarrow \mathbb{R}$ used in the CM is
continuous. Then for every $x \in X$ and $\epsilon > 0$ there exist $\delta > 0$
such that $f( B(x,\epsilon) ) \subset B(f(x),\delta)$. In other words, points covered by each ball in
BM are supposed to have similar values of the function $f$. In this case we may hope for some correspondence between the graphs of CM and BM. Before formalizing it, let us consider one dimensional example presented in the Figure~\ref{fig:difference_between_mappers}. It indicates that the
connectivity of the two graphs may be drastically different if the grow of the function is too large compared to
the density of sampling.
In those cases there may exist disconnected regions of CM that are connected in BM, as BM do not depend on the lens function, but only on proximity of points in $X$.
To avoid this situation let us assume that $f$ is uniformly continuous, i.e. the constants $\epsilon$ and $\delta$ are universal for every $x \in X$. Moreover let us assume that $\epsilon$ is the (initial) radius of balls in multi-scale BM construction. Further let us assume that the cover of $\mathbb{R}$ in CM consist of intervals of diameter $6\delta$ that overlap on the sub-intervals of the length $\delta$ and that the clustering algorithm used in CM is a single linkage clustering~\cite{single_linkage} with the parameter $\epsilon$. For comparison we will use a BM graph with a parameter $\epsilon$.

Firstly let us take $x,y \in X$ that in BM are covered by the balls $B(c_x,\epsilon)$ and $B(c_y,\epsilon)$ such that $B(c_x,\epsilon) \cap B(c_y,\epsilon) \neq \emptyset$\footnote{Meaning, there exist a point in $X$ that is covered by both the balls, and therefore there is an edge between $c_x$ and $c_y$ in the BM graph.} in the BM. We will show that in this case either $x$ and $y$ belongs to the same vertex in CM, or they belong to two vertices connected in CM. As $x \in B(c_x,\epsilon)$ and $y \in B(c_y,\epsilon)$, $d(x,y) < 4\epsilon$ that implies $d( f(x),f(y) ) < 4\delta$. Consequently $x$ and $y$ will be mapped by $f$ to one or at most two intersecting elements of the range cover.

Let us suppose first that they are mapped to the same element of range cover. Then as $x \in B(c_x,\epsilon)$, $d( x,c_x ) < \epsilon$. Analogously $d(y,c_y)<\epsilon$. Lastly as $B(c_x,\epsilon)\cap  B(c_y,\epsilon) \neq \emptyset$, $d(c_x,c_y) < \epsilon$. Consequently the single linkage clustering with a parameter $\epsilon$ will put $x$ to the same cluster as $y$, therefore they will be in the same vertex of CM.

Let us now assume that $x$ and $y$ are mapped to two intersecting elements $A$ and $B$ of
the range cover. If they are mapped to $A \cup B$, then the conclusions of the previous point apply as they will end up in the same vertex of $A$ and $B$. Therefore without loose of generality let us assume that $f(x) \in A \setminus B$ and $f(y) \in B \setminus A$.
As all distances $d(x,c_x),d(c_x,c_y)$ and $d(c_y,y)$ are bounded by $\epsilon$, at least one of
$c_x$ or $c_y$ have to be mapped to $A \cap B$. This vertex will be then connected (possibly via the
other one) both to $x$ and $y$. Therefore the vertices of the CM mapper graph $x$ and $y$ belong to
will be connected what finish the proof.

Under the assumptions listed above let us ask the reverse question: suppose we are given two vertices $x,y \in X$ which are grouped in one, or two nodes $A,B$ of CM graph, such that there is an edge between $A$ and $B$. Will they also be in neighboring nodes (elements of cover) of a BM graph?

Let us firstly assume that $x$ and $y$ belong to the same vertex $A$ of CM graph. In that case
they are mapped by the function $f$ to the same element $C$ of the range cover. Therefore
$d( f(x),f(y) ) < 6 \delta$, so their values are close. As they are both in the same node of CM
graph, then there exist a path composed of the edges of a length at most $\epsilon$ in $f^{-1}(C)$
joining $x$ and $y$. However the distance between them may be arbitrary large, we only know that both
$x$ ad $y$ are connected by according to a singe linkage clustering with a parameter $\epsilon$ in $f^{-1}(C)$.

The case when $x$ and $y$ belong to two neighboring
vertices of CM graph can be solver using the same argument as above, and he conclusion
holds. Therefore also in this case, the vertices in BM graph may be far away, but they will be connected.
Once the BM graph is colored using a function $f$, the user will be able to determine that both
$x$ and $y$ belong to the same connected component of the fiber of $f$.

To conclude, points that are close in the BM are expected to be close in the CM. Points close in CM will be connected in BM, but may be arbitrary far away. As a consequence, BM may bring more accurate information about the data, but at the same time, this information may be much more difficult to interpret. As a consequence the BM technique should be considered an additional technique that can be used together with CM for exploratory data analysis.
\subsection{Dimension estimation}
\label{sec:dimension-estimations}
\begin{figure}[ht]
\vskip 0.2in
\begin{center}
\includegraphics[scale=0.5]{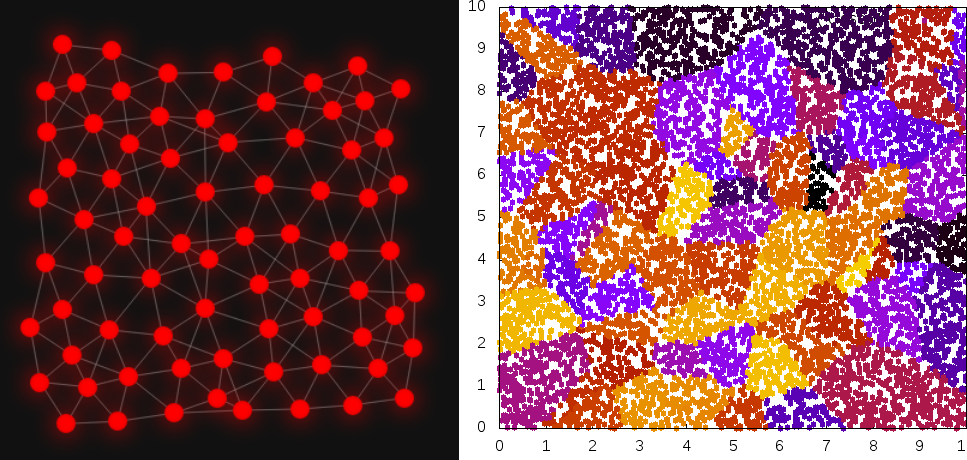}
\caption{\emph{(left) BM graph of $10000$ points randomly uniformly sampled from $[0,10]^2$ (presented on the right). (right) The initial point cloud colored by the closest ball (aka Voronoi cells restricted to the point cloud). The radius in the BM construction is $1$.}}
\label{fig:dmension_estimation_2_d_example_voronoi}
\end{center}
\vskip -0.2in
\end{figure}

As we have seen in Section~\ref{sec:relaton_to_standard_mapper}, BM graph is expected to contain more vertices than CM graph and may carry some additional information. As an example of such a scenario in this section we will analyze a number of neighbors of vertices in BM graphs. We will argue that by dong so one may get a proxy for the local dimension of the dataset. As a motivation let us consider a two dimensional data set presented in the Figure~\ref{fig:dmension_estimation_2_d_example_voronoi}

As one can read from the Figure~\ref{fig:dmension_estimation_2_d_example_voronoi} a typical vertex of the BM graph that is inside the region have six or seven neighbors. A vertex of the BM graph in the boundary (but not a corner) have four or five, and vertex in the corner three neighbors. Those numbers do not vary with radius as long as it is smaller than the diameter of the dataset and large enough so that defects resulting from non uniform sample are not visible.

\begin{figure}[ht]
\vskip 0.2in
\begin{center}
\includegraphics[scale=0.5]{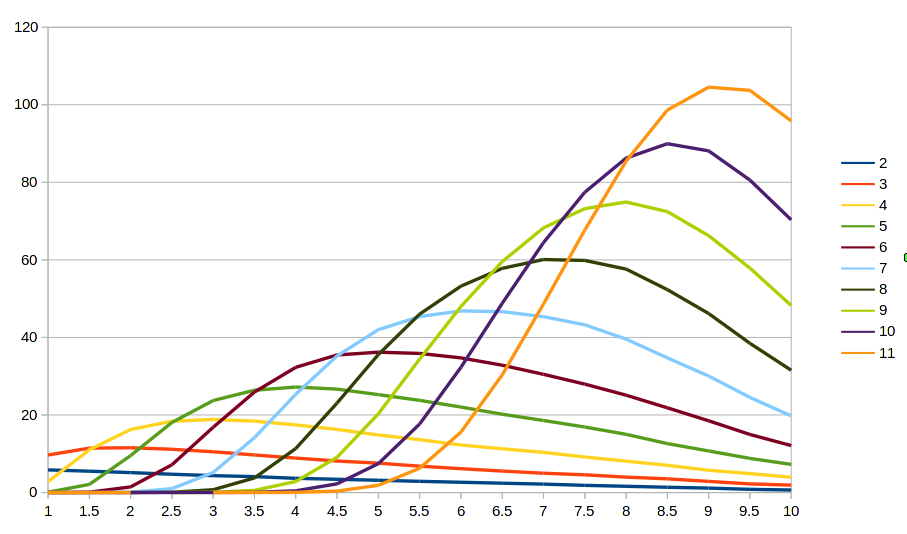}
\caption{Average number of neighbours per a node of BM graph obtained from $10000$ points sampled from $[0,10]^d$. $\epsilon$ is given in the $x$ axis.}
\label{fig:dmension_estimation_2_d_example}
\end{center}
\vskip -0.2in
\end{figure}

To numerically verify this hypothesis we have run $1000$ times  the following experiment and averaged results. We have sampled $10 000$ points uniformly iid from $[0,10]^d$, for $d \in \{2,\ldots,10\}$.
For each of the samples we have constructed a number of BM graphs for $r \in \{1,1.5,2,2.5,3,\ldots,9,9.5,10\}$. For each BM graph we have computed an average number of neighbours for all vertices. The results of the computations are presented in the
Figure~\ref{fig:dmension_estimation_2_d_example}. As one can observe, the plots are very flat for
low dimensions: most of the variation is placed for small and large values of $\epsilon$. That
indicate that excluding the case of a too small radii (where we capture local defects in the
sample) and too large radii (where the radii are comparable to the size of the point cloud when every vertex gets connected to any other vertex) the number of neighbors is roughly constant, and different for different dimensions.

The phenomena mentioned above is far less visible for higher dimensions where the flat part of curves get much narrower as the dimension grows. We believe that this phenomena is triggered by decreasing density of points $[0,10]^d$ as $d$ grows. In this case we can see artifacts of sampling for much longer. It is a manifestation of a curse of dimensionality. Note hewer that we still obtain well separated curves, and therefore an average number of neighbors may give an estimate of the local dimension of the dataset.

We would like to stress that the estimates we get by analyzing an average number of neighbors in BM graphs only give an upper bound for a dimension: for some selections of the ball's centers one can have a data sets of a low intrinsic dimension such that every vertex of its BM graph may have arbitrary many neighbors. Consequently the data set from a perspective of BM seems to have high intrinsic dimension. However, for a data set of a given dimension $d$, there is a certain minimal number of balls required to locally cover the neighborhood a ball in this data set. This will be the minimal number of neighbors of a vertex in a BM graph. This number is related to a box counting dimension of the dataset (at a certain resolution). Therefore it is not possible that a high dimensional dataset is classified as low dimensional one by analyzing the average number of neighbors in BM graph. Consequently, this technique can only be used to provide an upper bound for a dimension of a data set.
\subsection{Guarantees on connectivity of BM and its relation to persistent homology.}
\label{sec:Guarantees_on_connectivity}
In this section we provide some guarantees of the connectivity of the multi scale BM. Given a point cloud $X$ and a sequence of radii $\epsilon < \epsilon_1 < \ldots \epsilon_n$ let us assume that the multi scale BM have been constructed. Let $B$ denote the collection of ball's centers used in this construction. Then we have the following sequence of inclusions:
\[\bigcup_{b \in B}B(b,\epsilon) \subset \bigcup_{x \in X}B(x,\epsilon) \subset \bigcup_{b \in B}B(b,2\epsilon)\]
The first inclusion follows directly from the fact that $B \subset X$. For the second we use the
fact that $B \subset X$ is an $\epsilon$ net in $X$, i.e. for every $x \in X$ there exist $b \in B$
such that $d(x,b) < \epsilon$. Let us take any $y \in \bigcup_{x \in X} B(x,\epsilon)$. That implies
existence of $x \in X$ such that $d(x,y) < \epsilon$. Also, for the $x$ there exist $b \in B$
that $d(x,b) < \epsilon$. Summing up, $d(y,b) \leq d(y,x) + d(x,b) < 2\epsilon$. Consequently $y \in
\bigcup_{b \in B} B(b , 2\epsilon)$, what proves he second inclusion.

Consequently, building up on the multi scale BM, as described in Section~\ref{sec:mult-scale-approach}, if a feature is visible on the level $\epsilon$ and $2\epsilon$, then it is also present in the $\epsilon$ neighborhood of all the points in $X$.

One can also observe, using the same argument as above, that the sequence of inclusions can be extended to:
\begin{multline}
\nonumber
\bigcup_{b \in B}B(b,\epsilon) \subset \bigcup_{x \in X}B(x,\epsilon) \subset \bigcup_{b \in B}B(b,2\epsilon) \subset \bigcup_{x \in X}B(x,2\epsilon)
\subset \bigcup_{b \in B}B(b,3\epsilon) \subset \ldots
\end{multline}

That yields the interleaved sequence of homology groups:
\begin{center}
\includegraphics[scale=0.2]{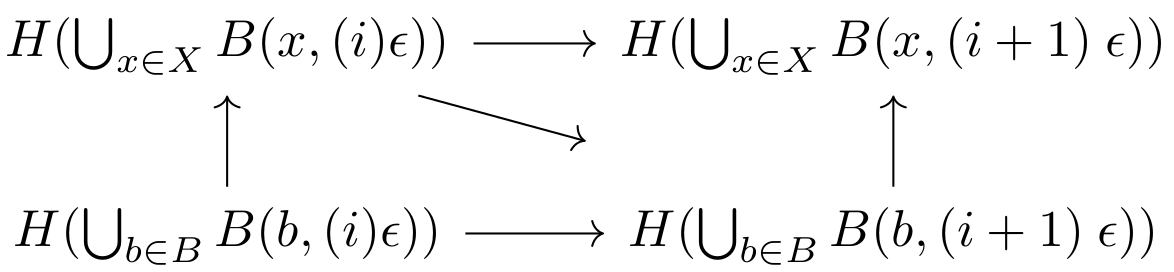}
\end{center}

Therefore once we compute the higher dimensional nerve of the cover, we can approximate the PH of the space $X$ from the nerve complex of points in $B$. That observation shows a relation between BM and PH, and shows how the sub-sampling procedure used in BM construction can be used to approximate PH.

\subsection{BM graph for noisy data}
In this section analyze sensitivity of BM graph construction with respect to Hausdorff-bounded and uniform noise.

To speak about the Hausdorff bounded noise, let us recall a Hausdorff metric. For two point clouds $X$ and $Y$ their Hausdorff distance $d_{H}(X,Y) = max\{ sup_{x \in X} inf_{y \in Y} d(x,y), sup_{y \in Y} inf_{x \in X} d(x,y)  \}$, where $d(x,y)$ denote a distance between points (typically Euclidean distance). The definition implies that if $d_{H}(X,Y) < \delta$, then for every point $x \in X$ there exist $y \in Y$ not farther away than $\delta$ from $x$. And the oher way around. Given this observation we have the natural sequence of inclusions for point clouds $X$ and $Y$. Suppose $B$ is an $\epsilon$ net in $X$.
\begin{multline}
\nonumber
\bigcup_{ b \in B } B(b,\epsilon) \subset \bigcup_{ x \in X } B(x,2\epsilon) \subset \bigcup_{ y \in Y } B(y,2\epsilon+\delta) \subset \subset \bigcup_{ x \in X } B(x,2\epsilon+2\delta)  \subset \bigcup_{ b \in B } B(b,3\epsilon+2\delta)
\end{multline}
Therefore if we have an estimate $\delta$ in Hausdorff metric on the noise, we can trust the features that persist between level $\epsilon$ and $3\epsilon+2\delta$ in BM graph.

\begin{figure}[ht]
\vskip 0.2in
\begin{center}
\includegraphics[scale=0.35]{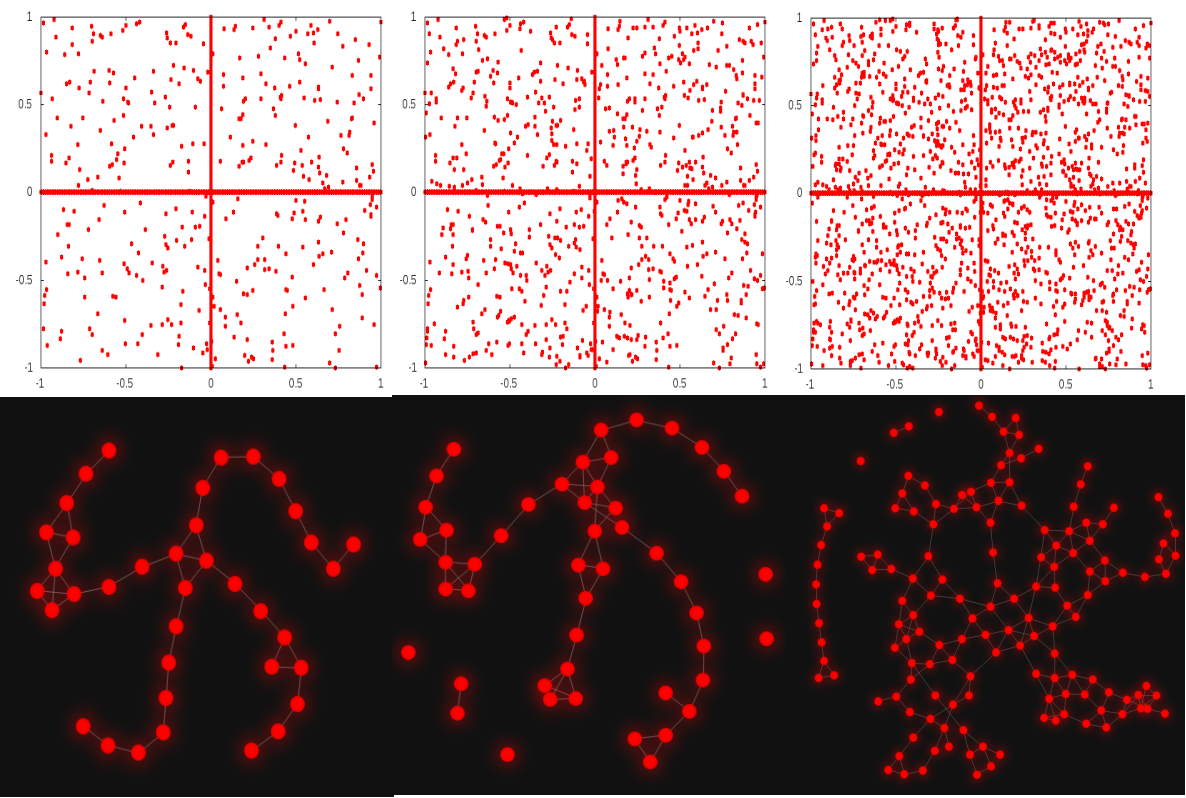}
\caption{\emph{A signal consist of the x-shape. Top left to right it is equipped with $50, 100, 150 \%$ of uniform noise. Bottom row indicate BM graphs in which cover elements that cover small number of points in $X$ have been removed.}}
\label{fig:noisy_x_junction}
\end{center}
\vskip -0.2in
\end{figure}

This type of noise is typical to a reconstruction or numerical error. Another type or nose, presented in the Figure~\ref{fig:noisy_x_junction}, is a uniform noise. Once we run the Algorithm~\ref{alg:ball_mapper_final} for point cloud with a uniform noise, we will obtain ball centers all over the support of signal and noise. However we can observe that the density of points in balls that belong to the signal is superior to the density of points to the balls in the noisy area. We can use this observation to clear out the regions of lower density. In the Figure~\ref{fig:noisy_x_junction} we are removing cover elements which cover less than a fixed number $n$ of points in $X$. Some strategies to estimating the number $n$ will be discussed in the full version of the paper.

\section{Experiments}
\begin{figure}[ht]
\vskip 0.2in
\begin{center}
\includegraphics[scale=0.5]{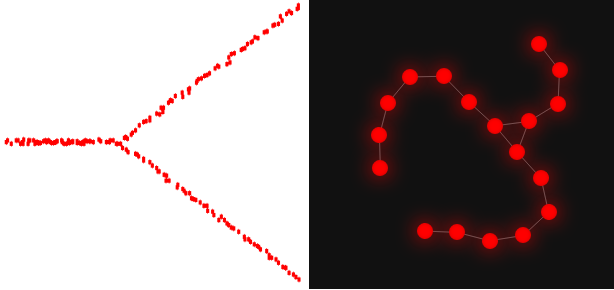}
\caption{\emph{Y-junction and its BM graph}}
\label{fig:Y_junction}
\end{center}
\vskip -0.2in
\end{figure}

In this section we will provide a number of examples of usage of the BM algorithm. Let us start with a standard two dimensional example of Y-junction presented in the Figure~\ref{fig:Y_junction}.
In this case the resulting graph will be very similar to CM graph for lens being a projection to x coordinate. Note however that if this data set is embedded (even linearly) in a high dimensional space, we will not be able to get the shape for every lens function. On the other hand, BM depends only on the geometric organization of points in the data set and will always recover the Y junction structure.

\begin{figure}[ht]
\vskip 0.2in
\begin{center}
\includegraphics[scale=1]{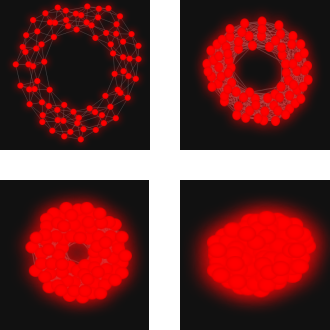}
\caption{\emph{BM construction on two dimensional torus for increasing sequence of radii.}}
\label{fig:tore_dataset}
\end{center}
\vskip -0.2in
\end{figure}

In the second example $1300$ points has been sampled from a 2 dimensional torus embedded to $\mathbb{R}^3$. A multi-scale BM construction has been run on it. The result can be seen in the Figure~\ref{fig:tore_dataset}.

Over here we can already observe a substantial differences between the CM and BM. While for the height function lens the resulting CM graph will consist of a single dominant loop with potentially a few out-layers, in the BM we clearly see the volumetric structure of data. For the BM graph the averaged number of vertices neighbors oscillates in between 5-6 which gives and indication that the set is locally two dimensional.

Let us point out that as the construction of BM depends on distances between points, the results presented above will no change if the data are embedded in high dimensional space as long as the distances are not perturbed too much. This is typically not the case for CM.

\begin{figure}[ht]
\vskip 0.2in
\begin{center}
\includegraphics[scale=0.7]{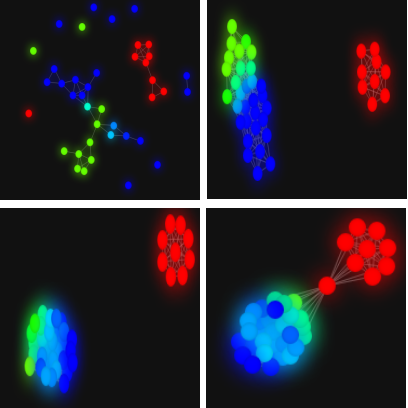}
\caption{\emph{Mult scale BM construction on Iris dataset for initial radius $0.5$ and the radii $0.9,1.6,1.8$.}}
\label{fig:iris_dataset}
\end{center}
\vskip -0.2in
\end{figure}

Let us now analyze a few standard data sets, starting from Iris data set obtained from UCI Machine Learning Repository~\cite{iris_UCI} originated from~\cite{iris_original}. It describe the taxonomy of three different types of iris plant. Please consult Figure~\ref{fig:iris_dataset} for the result. On the upper left, the result for $\epsilon = 0.5$. Over there we start seeing clusters forming, and for $\epsilon = 0.9$ on the upper right they are fully formed. We can see that the class Iris-setosa (in red) is clearly separated from the class Iris-versicolor (green) and Iris-virginica (blue), while the green and blue classes continuously transfer to each other. The separation between red on one side, and green and blue on the other persist for long, through $\epsilon = 1.6$ on the bottom left when we see green and blue classes fusing together. They get merge for $\epsilon = 1.8$ on the bottom right. We can see that the red class gets merged to the main cluster via the green class.

\begin{figure}[ht]
\vskip 0.2in
\begin{center}
\includegraphics[scale=0.5]{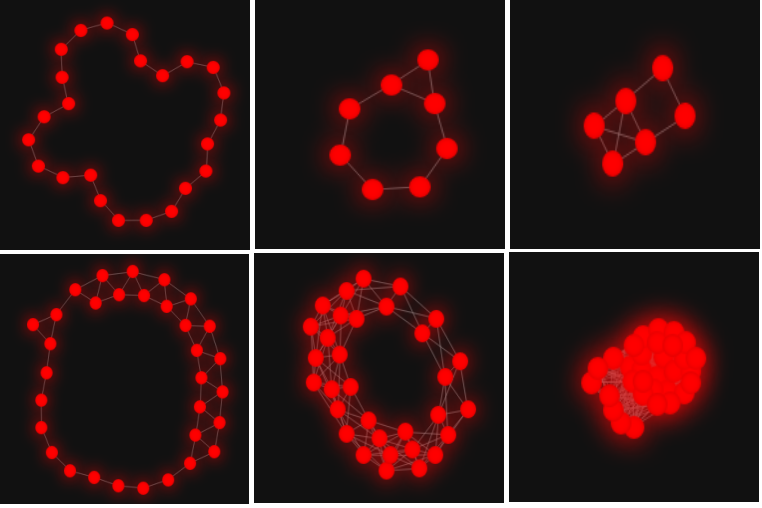}
\caption{\emph{In the top row we have constructed the BM graphs independently for $\epsilon=5000, 8000$ and $9000$. The cycle we can see over there is created just before $\epsilon = 5000$. Bottom line shows the multi-scale BM where the centers of balls has been selected for $\epsilon = 4000$. The picture shows the multi-scale BM graph for $\epsilon = 5000, 8000$ and $9000$.}}
\label{fig:lucky_cat}
\end{center}
\vskip -0.2in
\end{figure}

Our next example involve a data set where we can expect to see one dimensional structure. It is a \emph{lucky cat} data set taken from~\cite{columbia_image_dataset}. It consist of seventy two $128 \times 128$ images of a sculpture of a cat taken from a different angles as described in the Figure 1 in~\cite{lucky_cat_exmplanation}. Given this we can expect that the data form a circle in $128 \times 128 = 16384$ dimensional space. Indeed, as we can observe in the Figure~\ref{fig:lucky_cat} this is the case. As we can observe, the top row give us an idea about the size of the cycle in the ambient space. The bottom row gives an idea about some non uniformity in sampling. Note that there may not be a cycle in the lower right picture, as there are many more centers of balls over there, and they may already cover the internal hole.

\section{Implementation}
A prototype $C++$ implementation of the presented method that uses the Kepper Mapper~\cite{keppler_mapper} for visualization is available upon request. Once tested it will be distributed publicly. 

\section{Final remarks}
The presented technique generalize in an obvious way to \emph{weighted graphs}. The input is not restricted to point cloud. We can use any discrete metric space with a distance matrix or similarity measure. All those extensions will be discussed in a full version of the paper.

\section*{Acknowledgements}
I want to thank Gunnar Carlsson and Tan Li for helpful discussion and feedback.

\bibliographystyle{icml2018}

\end{document}